\documentclass[12pt]{amsart}
\usepackage{lipsum}
\usepackage{amssymb}
\usepackage{enumerate}
\usepackage{amsthm}
\usepackage{tikz-cd}
\usepackage{amsmath}
\usepackage[mathscr]{euscript}
\usepackage[utf8]{inputenc}
\usepackage[english]{babel}
\usepackage{hyperref}
\hypersetup{
    colorlinks=true,
    linkcolor=blue,
    citecolor=red,
    filecolor=red,      
    urlcolor=violet,
}

\makeatletter
\@namedef{subjclassname@2020}{%
  \textup{2020} Mathematics Subject Classification}
\makeatother

\newtheorem{theorem}{Theorem}[section]

\theoremstyle{definition}

\newtheorem{exm}[theorem]{Example}

\newtheorem{conj}[theorem]{Conjecture}

\theoremstyle{remark}
\newtheorem{xrem}[theorem]{Remark}
\newtheorem{prop}[theorem]{Proposition}
\newtheorem{defi}[theorem]{Definition}

\newtheorem{corl}[theorem]{Corollary}

\numberwithin{equation}{section}

\DeclareMathOperator{\Nef}{{Nef}}

\DeclareMathOperator{\rk}{{rk}}

\DeclareMathOperator{\Amp}{{Amp}}

\DeclareMathOperator{\Sym}{{Sym}}

\DeclareMathOperator{\End}{{\mathcal{E}nd}}
\DeclareMathOperator{\Div}{{Div}}
\DeclareMathOperator{\Proj}{{Proj}}

\DeclareMathOperator{\NE}{{NE}}

\begin{document}
\baselineskip=15pt

\title[On Serrano's conjecture on Projective bundles]{On Serrano's conjecture on Projective bundles}
\subjclass[2020]{Primary 14J40, 14J45, 14J60}
\keywords{Strictly nef, Projective bundle, Serrano's conjecture}

\author{Snehajit Misra}

\address{Harish Chandra Research Institute, India.}
\email[Snehajit Misra]{misra08@gmail.com}

\begin{abstract}
Serrano’s conjecture asserts that if $D$ is a strictly nef divisor on a projective variety $X$ of dimension $n$, then $K_X + tD$ is ample for $t > n + 1$. In this article, we investigate this conjecture for  strictly nef divisors on projective bundles over higher dimensional smooth projective varieties.
\end{abstract}

\maketitle

\section{Introduction}
A line bundle $L$ on a projective variety  $X$ is called \it strictly nef \rm if $L\cdot C > 0$  for every irreducible curve $C$ in $X$. A  divisor  $D$ on $X$ is called strictly nef if the corresponding line bundle $\mathcal{O}_X(D)$ is strictly nef.
The Nakai-Moishezon-Kleimann criterion asserts that a line bundle $L$
is ample if and only if $L^{\dim V}\cdot V >0$ for every positive dimensional irreducible subvarieties $V$ in $X$. Hence any ample line bundle is strictly nef. However Mumford constructed a number of  strictly nef but non-ample line bundles over ruled surfaces (e.g. \cite{H70}).  These strictly nef line bundles are tautological line bundles of stable vector bundles of degree zero over smooth curves of genus $g \geq  2$.

 Campana and Peternell proposed in \cite[Problem 11.4]{CP91} the following conjecture, which is still unsolved.
\begin{conj}\label{conj1}
Let $X$ be a projective manifold. If $-K_X$ is strictly nef, then $X$ is Fano.
\end{conj}
The dual version of Conjecture \ref{conj1} follows from  the abundance conjecture.
\begin{conj}
Let $X$ be a projective manifold. If $K_X$ is strictly nef, then $K_X$ is ample.
\end{conj}
Analogous to the Fujita conjecture, Serrano asserts the following
conjecture, which is a generalization of Conjecture \ref{conj1}.
\begin{conj}
 If $D$ is strictly nef on a projective variety $X$ of dimension $n$ then $K_X+tD$ is ample for $t> n+1$.
\end{conj}
This conjecture has been solved by Serrano for projective surfaces in \cite{S95}. For the progress on projective threefolds and higher dimensional manifolds,  see \cite[Theorem 0.4]{CCP08}. More specifically, the following is proved in  \cite[Theorem 0.4]{CCP08}.
\begin{theorem}
Let $X$ be a projective manifold of dimension $n$ and $L$ a strictly nef line bundle on $X$. Then $K_X + tL$ is ample if $t > n + 1$ in the following cases.
\begin{enumerate}
    \item  $\dim X = 3$ unless (possibly) $X$ is Calabi-Yau with $L\cdot c_2 = 0$.
    \item $\kappa(X) \geq n-2$
    \item  $\dim \alpha(X) \geq n - 2$, with $\alpha : X \longrightarrow A$ the Albanese map.
\end{enumerate}
\end{theorem}

%By Mori's cone theorem, Serrano's conjecture is equivalent to the condition that $K_X^{\perp}\cap D^{\perp} = 0$ in Mori closed cone $\overline{\NE}(X)$ of curves in $X$.

A vector bundle $E$ on a smooth projective variety $X$ is called striclty nef (resp. ample) in the sense of Hartshorne if the tautological line bundle $\mathcal{O}_{\mathbb{P}(E)}(1)$ on $\mathbb{P}(E)$ is striclty nef (resp. ample). Hence by Nakai-Moishezon-Kleimann criterion, any ample vector bundle is striclty nef. In this article, we investigate this conjecture for strictly nef divisors on projective bundles over higher dimensional smooth projective varieties.

 Let $E$ be a strictly nef vector bundle of rank $r$ on a smooth projective irreducible variety $X$ of dimension $n$, and $\pi : \mathbb{P}(E)\longrightarrow X$ be the projection. Then consider the divisor $D\equiv m\xi+\pi^*N$ for any positive integer $m>0$ and for any nef divisor class $N$ on $X$.
  Note that $\pi^*N$ is nef and the numerical class $\xi$ of the tautological bundle $\mathcal{O}_{\mathbb{P}(E)}(1)$ is strictly nef. Hence for any irreducible curve $C\subseteq X$, we have $$D\cdot C = \bigl(m\xi+\pi^*D\bigr)\cdot C > 0.$$
  This shows that $D$ is strictly nef. Throughout this article, we will prove Serrano's conjecture for such strictly nef divisors, i.e. we will show that $K_{\mathbb{P}(E)}+tD$ is ample for $t\geq \dim\bigl(\mathbb{P}(E)\bigr)+2 = n+r+1$ for strictly nef divisors $D$ as mentioned above. In particular, we prove the following results.

  \begin{theorem}
Let $E$ be a strictly nef vector bundle of rank $r$ on a complex projective curve $C$, and $\pi : X = \mathbb{P}(E)\longrightarrow C$ be the projection.  Consider the strictly nef divisor $D=m\xi+\pi^*N$ on $X$ for some positive integer $m>0$ and for some nef divisor $N$.
 Then Serrano's conjecture is true for the strictly nef divisor $D$.
\end{theorem}
\begin{theorem}
Let $\pi : \mathbb{P}(E)\longrightarrow X$ be a projective bundle corresponding to a strictly nef  vector bundle $E$ of rank $r\geq 2$ on a smooth projective variety $X$ of dimension $n$. Further assume that one of the following holds :
\begin{enumerate}
   \item $K_X+\det(E)$ is strictly nef.
  \item $E$ is ample and $K_X+\det(E) \equiv 0.$
\end{enumerate}
  Let $D\equiv m\xi+\pi^*N$ be a strictly nef divisor  on $\mathbb{P}(E)$ for some $m>0$ and some nef divisor class $N$ on $X$. Then Serrano's conjecture is true for the strictly nef divisor $D$ on $\mathbb{P}(E)$.
  \end{theorem}
\begin{theorem}
Let $E$ be a Fano ample vector bundle of rank 2 on $\mathbb{P}^n$, ( i.e. $\mathbb{P}(E)$ is Fano) and $\pi : \mathbb{P}(E)\longrightarrow \mathbb{P}^n$ be the projection. Let $D \equiv m\xi + \pi^*N$ be a strictly nef divisor on $\mathbb{P}(E)$ for some $m > 0$ and some nef divisor class $N$ on $\mathbb{P}^n$.
Then Serrano’s conjecture is true for the strictly nef divisor $D$ on $\mathbb{P}(E)$.
\end{theorem}

\section{Preliminaries}
 Throughout this article, all the algebraic varieties are assumed to be irreducible and defined over the field of characteristic zero. The set of non-negative  real numbers (resp. non-negative rational numbers) will be denoted by $\mathbb{R}_{\geq 0}$ (resp. $\mathbb{Q}_{\geq 0}$). The cannonical line bundle of a projective variety $X$ will be denoted by $K_X$.

 In this section we recall basic definition of closed cones of curves and nef cones of divisors . We also recall the notion of semistability of vector bundles on smooth projective varieties.
 \subsection{Nef cone of divisors and closed cone of curves}
 Let $X$ be a smooth projective variety of dimension $n$ and $\mathcal{Z}_1(X)$ (respectively $\mathcal{Z}^1(X) = \Div(X)$) denotes the free abelian group generated by $1$-dimensional (respectively $1$-co dimensional) closed subvarieties on $X$.
 Two cycles $Z_1$,$Z_2 \in \mathcal{Z}_1(X)$ are said to be numerically equivalent, denoted by $Z_1\equiv Z_2$ if $Z_1\cdot \gamma =  Z_2\cdot\gamma $
  for all $\gamma \in \mathcal{Z}^1(X)$. The \it numerical groups \rm  $ N_1(X)_{\mathbb{R}}$ are defined as the quotient of
  $\mathcal{Z}_1(X) \otimes \mathbb{R}$ modulo numerical equivalence. Let
  \begin{center}
  $\Div^0(X) := \bigl\{ D \in \Div(X) \mid D\cdot C = 0 $ for all curves $C$ in $X \bigr\} \subseteq \Div(X)$.
  \end{center}
 be the subgroup of $\Div(X)$ consisting of numerically trivial divisors. The quotient\\ $\Div(X)/\Div^0(X)$ is called the N\'{e}ron Severi group of $X$, and is denoted by $N^1(X)_{\mathbb{Z}}$.
   The N\'{e}ron Severi group  $N^1(X)_{\mathbb{Z}}$ is a free abelian group of finite rank.
 Its rank, denoted by $\rho(X)$ is called the Picard number of $X$.
 In particular,  $N^1(X)_{\mathbb{R}}  := N^1(X)_{\mathbb{Z}} \otimes \mathbb{R}$ is a finite dimensional real vector space of dimension $\rho(X)$. The real vector space $N^1(X)_{\mathbb{R}}$ is called the real N\'{e}ron
 Severi group of $X$.  The numerical class of a divisor $D$ will be denoted by  $\bigl[D\bigr]$. Two numerically equivalent divisors  $D_1$ and $D_2$ will be denoted by $D_1\equiv D_2$.
The intersection product induces a perfect pairing
 \begin{align*}
N^1(X)_{\mathbb{R}} \times N_1(X)_{\mathbb{R}} \longrightarrow \mathbb{R},
\end{align*}
which implies $N^1(X)_{\mathbb{R}} \cong (N_1(X)_{\mathbb{R}})^\vee$. The convex cone generated by the set of all effective $1$-cycles in $N_1(X)_\mathbb{R}$ is denoted by $\NE(X)$ and its closure $\overline{\NE}(X)$ is called the \it closed cone of curves \rm in $X$. The \it nef cone \rm of divisors are defined as follows :
\begin{align*}
 \Nef(X) := \bigl\{ \alpha \in N^1(X)_{\mathbb{R}} \mid \alpha \cdot \beta \geq 0 \hspace{2mm} \forall \beta \in \overline{\NE}(X)\bigr\}.
\end{align*}
The interior of $\Nef(X)$  is called ample cone of  divisors in $X$, and is denoted by $\Amp(X)$. By \cite[Theorem 1.4.29]{L1}, an $\mathbb{R}$-divisor class $\alpha\in N^1(X)_{\mathbb{R}}$ is ample if and only if $\alpha\cdot \gamma > 0$ for all $\gamma \in \overline{\NE}(X)  \setminus\{0\}$.

\subsection{Semistability of Vector bundles}
Let $X$ be a smooth complex projective variety of dimension $n$ with a fixed ample divisor $H$ on it.
For a torsion-free coherent sheaf $\mathcal{G}$ of rank $r$ on $X$, the $H$-slope of $\mathcal{G}$ is defined as
\begin{align*}
\mu_H(\mathcal{G}) := \frac{c_1(\mathcal{G})\cdot H^{n-1}}{r} \in \mathbb{Q}.
\end{align*}
A torsion-free coherent sheaf $\mathcal{G}$ on $X$ is said to be  $\mu_H$-semistable if $\mu_H(\mathcal{F}) \leq \mu_H(\mathcal{G})$ for
all subsheaves $\mathcal{F}$ of $\mathcal{G}$.
A vector bundle $E$ on $X$ is called $\mu_H$-unstable if it is not  $\mu_H$-semistable. For every vector bundle $E$ on $X$, there is a unique filtration
\begin{align*}
 0 = E_d \subsetneq E_{d-1} \subsetneq E_{d-2} \subsetneq\cdots\subsetneq E_{1} \subsetneq E_0 = E
\end{align*}
of subsheaves of $E$, called the Harder-Narasimhan filtration of $E$, such that each successive quotient $E_i/E_{i+1}$ is $\mu_H$-semistable torsion-free sheaf for each $i \in \{ 0,1,2,\cdots,d-1\}$
and $$\mu_H\bigl(E_{d-1}/E_{d}\bigr) > \mu_H\bigl(E_{d-2}/E_{d-1}\bigr) >\cdots> \mu_H\bigl(E_{0}/E_{1}\bigr).$$
We define $Q_1:= E_{0}/E_{1}$ and $\mu_{\min}(E) := \mu_H(Q_1) = \mu_H\bigl(E_{0}/E_{1}\bigr)$ and $\mu_{\max}(E) := \mu_H\bigl(E_{d-1}/E_d\bigr).$ We always have the following inequality : $$\mu_{\min}(E) \leq \mu(E) \leq \mu_{\max}(E)$$ and $E$ is $\mu_H$-semistable if and only if $\mu_{\min}(E) = \mu_H(E) = \mu_{\max}(E)$.
\vspace{2mm}

\subsection{Nef cone of projective bundles} The projective bundle $\mathbb{P}_X(E)$ associated to a vector bundle
$E$ over a projective variety $X$ is defined as $\mathbb{P}_X(E) :=  \Proj\bigl(\bigoplus\limits_{m\geq  0} \Sym^m(E)\bigr)$ together with the projection map $$\pi : \mathbb{P}_X(E)\longrightarrow X.$$
We will simply write $\mathbb{P}(E)$ whenever the base space $X$ is clear from the context. The numerical class of the tautological line bundle $\mathcal{O}_{\mathbb{P}(E)}(1)$ will be denoted by $\xi$ unless otherwise specified.
\vspace{2mm}

Let $\pi:\mathbb{P}(E)\longrightarrow C$ be a projective bundle  on a smooth projective curve $C$. Fix the following notations :
$\xi\equiv \mathcal{O}_{\mathbb{P}(E)}(1)$ and $f\equiv$ a fiber of $\pi$. Then
\begin{align*}
 \Nef\bigl(\mathbb{P}(E)\bigr) = \Bigl\{a\bigl(\xi-\mu_{\min}(E)f\bigr)+yf\mid x,y\in \mathbb{R}_{\geq 0}\Bigr\}.
\end{align*}
Thus a vector bundle $E$ on a smooth curve $C$ is nef (resp. ample) if and only if $\mu_{\min}(E)\geq 0 $ (resp. $\mu_{\min}(E) >0$).
\vspace{2mm}

Next we recall the description of nef cone of projective bundles over higher dimensional from \cite{MR22}.
Let $E$ be a vector bundle of rank $\geq 2$ on a smooth complex projective variety $X$, and $\pi:\mathbb{P}(E)\longrightarrow X$ be the projection. For an irreducible curve $C$ in $X$, let $\eta_c:\tilde{C}\longrightarrow C$ be it's normalization. Consider the following fibre product diagram:
\begin{center}
 \begin{tikzcd}
\tilde{C}\times_X\mathbb{P}(E)=\mathbb{P}\bigl(\eta_c^*(E\vert_C)\bigr) \arrow[r, "\tilde{\eta_c}"] \arrow[d, "\tilde{\pi_c}"] & C\times_X \mathbb{P}(E) = \mathbb{P}(E\vert_C) \arrow[r, "j"] \arrow[d,"\pi_c"] & \mathbb{P}(E) \arrow[d,"\pi"] \\
\tilde{C} \arrow[r, "\eta_c"]          & C \arrow[r, "i"]                                         & X
\end{tikzcd}
\end{center}
Here $i$ and $j$ are inclusions, and we define $\psi_c:= j \circ \tilde{\eta_c}$.
We fix the following notations $\xi_{\tilde{c}}\equiv \mathcal{O}_{\mathbb{P}(\eta_c^*(E\vert_C))}(1)$ and $f_{\tilde{c}}\equiv$ a fiber of $\tilde{\pi}_c$, $\xi\equiv \mathcal{O}_{\mathbb{P}(E)}(1)$.
 %Then $N^1\bigl(\mathbb{P}(\eta_c^*(E\vert_C)\bigr)_{\mathbb{R}}$ is generated by $\xi_{\tilde{c}}$ and $f_{\tilde{c}}$.
 The map $\psi_c$ induces the map
 \begin{center}
 $\psi_{c}^{*}:N^1\bigl(\mathbb{P}(E)\bigr)_{\mathbb{R}}\longrightarrow N^1\bigl(\mathbb{P}(\eta_c^*(E\vert_C)\bigr)_{\mathbb{R}}$
 \end{center}
  such that
 \begin{center}
 $\psi_{c}^{*}(\xi)=\xi_{\tilde{c}}$ and $\psi_{c}^{*}\bigl(\pi^*L\bigr)= (L\cdot C)f_{\tilde{c}},$ for any $L\in N^1(X)_{\mathbb{R}}.$
 \end{center}
 %This shows that the map $\psi_{c}^{*}$ is surjective.
 As a consequence, we get the push forward map
 $$(\psi_c)_* : N_1\bigl(\mathbb{P}(\eta_c^*(E\vert_C)\bigr)_{\mathbb{R}} \longrightarrow N_1\bigl(\mathbb{P}(E)\bigr)_{\mathbb{R}}.$$
 Hence
 \begin{align*}
 \overline{\sum\limits_{C\in \Gamma}(\psi_{c})_{*} \Bigl(\overline{\NE}\bigl(\mathbb{P}(\eta_c^*(E\vert_C))\bigr)\Bigr)} \subseteq \overline{\NE}\bigl(\mathbb{P}(E)\bigr),
\end{align*}
where the sum is taken over the set $\Gamma$ of all irreducible curves in $X$.

Next we consider the  numerical equivalence class $\bigl[\bar{C}\bigr]\in \overline{\NE}\bigl(\mathbb{P}(E)\bigr)$ of an irreducible curve $\bar{C}$ in $\mathbb{P}(E)$ which is not contained in any fibre of $\pi$.
Denote $\pi(\bar{C})=C$. Then, $\bar{C}\subseteq \mathbb{P}(E\vert_C)$. %Also, by Proposition \ref{prop3.2},
Then there exists a unique irreducible curve $C'\subseteq \mathbb{P}\bigl(\eta_c^*(E\vert_C)\bigr)$  such that $\tilde{\eta_c}(C') = \bar{C}$ and $\bigl(\psi_{c}\bigr)_*\bigl([C']\bigr)=\bigl[\bar{C}\bigr]$. Also, the numerical equivalence classes of curves in a fibre of $\tilde{\pi_{c}}$ maps to the numerical classes of curves in a fibre of $\pi$ by
 $\bigl(\psi_{c}\bigr)_*$.
Hence, we have
\begin{align*}
 \overline{\NE}\bigl(\mathbb{P}(E)\bigr) = \overline{\sum\limits_{C\in \Gamma}\bigl(\psi_{c}\bigr)_* \Bigl(\overline{\NE}\bigl(\mathbb{P}(\eta_c^*(E\vert_C))\bigr)\Bigr)}
\end{align*}
 where $\Gamma$ is the set of all irreducible curves in $X$.
\vspace{2mm}

  We use the notations as above from now on. In particular, we will use the notation $\eta_c$ to denote the normalization map $\eta_c:\tilde{C}\longrightarrow C$ for any irreducible curve $C\subseteq X$. We now find the ample cone $\Amp\bigl(\mathbb{P}(E)\bigr)$ by applying duality to $\overline{\NE}\bigl(\mathbb{P}(E)\bigr)\setminus \{0\}$.

\begin{corl}\label{corl4.5}
Let $\pi : \mathbb{P}(E)\longrightarrow X$ be a projective bundle corresponding to a vector bundle of rank $r\geq 2$ on a smooth irreducible complex projective variety $X$.
\vspace{2mm}

Then $y_0\xi+\pi^*\gamma$ is in the ample cone $\Amp\bigl(\mathbb{P}(E)\bigr)$ if and only if
\begin{center}
$y_0>0$ and $y_0\mu_{\min}\bigl(\eta_c^*(E\vert_C)\bigr)+\gamma\cdot C >0$ for every irreducible curve $C\subseteq X$.
\end{center}
\begin{proof}
For an irreducible curve $C\subseteq X$, we recall that
$$\Nef\bigl(\mathbb{P}(\eta_c^*(E\vert_C))\bigr) = \Bigl\{a\bigl(\xi_{\tilde{c}}-\mu_{\min}(\eta_c^*(E\vert_C)f_{\tilde{c}}\bigr)+bf_{\tilde{c}}\mid a,b \in \mathbb{R}_{\geq 0}\Bigr\}.$$
Let $l_c := \deg(\eta_c^*(E\vert_C))-\mu_{\min}(\eta_c^*(E\vert_C)).$ Applying duality we get
$$\overline{\NE}\bigl(\mathbb{P}(\eta_c^*(E\vert_C))\bigr) = \Bigl\{a\xi^{r-2}_{\tilde{c}}f_{\tilde{c}}+b(\xi^{r-1}_{\tilde{c}}-l_c\xi^{r-2}_{\tilde{c}}f_{\tilde{c}})\mid a,b\in \mathbb{R}_{\geq 0}\Bigr\}.$$
Thus we have $\overline{\NE}\bigl(\mathbb{P}(E)\bigr)$ is generated by
\begin{center}
  $\Bigl\{ \xi^{r-2}F$ , $\xi^{r-1}\pi^*[C] -l_c\xi^{r-2}F \mid C\subseteq X$ is an irreducible curve in $X \Bigr\}$,
\end{center}
where $F$ is the numerical equivalence class of the fibre of the projection map $\pi : \mathbb{P}(E)\longrightarrow X$. Again applying duality we conclude that $y_0\xi+\pi^*\gamma \in \Amp\bigl(\mathbb{P}(E)\bigr)$ if and only if
\begin{center}
$y_0>0$ and $y_0\mu_{\min}\bigl(\eta_c^*(E\vert_C)\bigr)+\gamma\cdot C >0$ for every irreducible curve $C\subseteq X$.
\end{center}
\end{proof}
\end{corl}
\begin{xrem}\label{xrem1}
 \rm By applying duality to $\overline{\NE}\bigl(\mathbb{P}(E)\bigr)$, we can show similarly that $y_0\xi+\pi^*\gamma \in \Nef\bigl(\mathbb{P}(E)\bigr)$ if and only if
\begin{center}
$y_0\geq 0$ and $y_0\mu_{\min}\bigl(\eta_c^*(E\vert_C)\bigr)+\gamma\cdot C \geq 0$ for every irreducible curve $C\subseteq X$.
\end{center}
\end{xrem}
\section{Main Results}
\begin{theorem}\label{thm4.6}
 Let $\pi : \mathbb{P}(E)\longrightarrow X$ be a projective bundle corresponding to a strictly nef  vector bundle $E$ of rank $r\geq 2$ on a smooth projective variety $X$ of dimension $n$. Further assume that Serrano's conjecture is true for any strictly nef divisor $D'$ on $X$.
 Then Serrano's conjecture is true for a strictly nef divisor $D \equiv m\xi+\pi^*D'$, where $D'$ is strictly nef on $X$, and $m>0$ is a positive integer.
 \begin{proof}
Note that $$K_{\mathbb{P}(E)} = -r\xi + \pi^*(K_X+\det(E)).$$
Thus $$K_{\mathbb{P}(E)} + tD \equiv \bigl(tm-r\bigr)\xi + \pi^*\bigl(K_X+tD'+\det(E)\bigr).$$
As Serrano's conjecture is true for the strictly nef divisor $D'$ on $X$, we have $K_X+tD'$ is ample for $t>n+1$.
As $D$ is strictly nef, we conclude $m>0$. Thus $tm-r > 0$ for $t> n+r$. By Corollary \ref{corl4.5}, it is enough to show that $$\bigl(tm-r)\mu_{\min}(\eta_c^*(E\vert_{C}))+\bigl(K_X+tD'+\det(E)\bigr)\cdot C >0$$ for every irreducible curve $C\subseteq X$.
As $E$ is strictly nef, both $E$ and $\det(E)$ are nef. Hence $\mu_{\min}(\eta_c^*(E\vert_{C}) \geq  0$ and $\det(E)\cdot C \geq 0$ for every irreducible curve $C\subseteq X$.

Now for any irreducible curve $C\subseteq X$, we observe that for $t>n+r>n+1$
\begin{align*}
 &\bigl(tm-r\bigr)\mu_{\min}(\eta_c^*(E\vert_{C}))+\bigl(K_X+tD'+\det(E)\bigr)\cdot C\\
 & = \bigl(tm-r\bigr)
 \mu_{\min}(\eta_c^*(E\vert_{C}))+\bigl(K_X+tD'+\det(E)\bigr)\cdot C\\
 & = \bigl(tm-r\bigr)\mu_{\min}(\eta_c^*(E\vert_{C})) +(K_X+tD')\cdot C + \det(E)\cdot C > 0.
\end{align*}
This completes the proof.
\end{proof}
\end{theorem}
\begin{theorem}
Let $E$ be a strictly nef vector bundle of rank $r$ on a complex projective curve $C$, and $\pi : X = \mathbb{P}(E)\longrightarrow C$ be the projection.  Consider the strictly nef divisor $D=m\xi+\pi^*N$ on $X$ for some positive integer $m>0$ and for some nef divisor $N$.
 Then Serrano's conjecture is true for the strictly nef divisor $D$.
\begin{proof}
We observe

$K_{\mathbb{P}(E)} + tD = \bigl(tm-r\bigr)\xi+\pi^*(K_C+\det(E)+tN)
 \vspace{1mm}$

 $= \bigl(tm-r\bigr)\bigl(\xi-\mu_{\min}(E)f\bigr)
 \vspace{1mm}$

 $\hspace{2mm} + \bigl(\deg(K_C)+\deg(E)-r\mu_{\min}(E) + tm\mu_{\min}(E)+t\deg(N)\bigr)f,
\vspace{2mm}$

where $f\equiv$ a fiber of $\pi$.
Note that for $t\geq \dim(\mathbb{P}(E))+2 = r+2$, we always have $\bigl(tm-r\bigr) > 0$. Hence  it is enough to show that $$\bigl(\deg(K_C)+\deg(E)-r\mu_{\min}(E) + tm\mu_{\min}(E)+t\deg(N)\bigr)>0$$ for $t\geq r+2$. As $E$ is strictly nef, in particular, $E$ is  nef. Thus $\mu_{\min}(E) \geq 0$ always. Also $N$ being nef, we have $\deg(N)\geq 0$. We always have $\deg(E)-r\mu_{\min}(E) \geq 0$.
\vspace{2mm}

 Now we consider the following three cases :
 \begin{itemize}
  \item Case (1) : Let the genus of the curve $C$ be at least 2. Then  we have $\deg(K_C) = 2g-2 >0$ and $K_C$ is ample. Thus
$K_{\mathbb{P}(E)}+tD$ is ample  in this case for $t\geq r+2.$
\item Case (2) : Let $C$ be an elliptic curve. Thus $\deg(K_C) =2g-2 = 0$. As $E$ is strictly nef vector bundle on an elliptic curve, by \cite[Theorem 3.1]{LOY19} we have $E$ is ample. Hence $\mu_{\min}(E)>0$ in this case.
   Thus for $t\geq r+3$
  $$\bigl(\deg(K_C)+\deg(E)-r\mu_{\min}(E) + tm\mu_{\min}(E)+t\deg(N)\bigr) >0.$$
\item Case (3) : Let $C=\mathbb{P}^1$ so that $\deg(K_C) = 2g-2 = -2$. As $E$ is strictly nef on $\mathbb{P}^1$, we have $E$ is ample. Thus $\mu_{\min}(E) > 0$.
Thus for $t\geq r+2 \geq 3$
  $$\bigl(\deg(K_C)+\deg(E)-r\mu_{\min}(E) + tm\mu_{\min}(E)+t\deg(N)\bigr) \geq t-2 > 0.$$
  \end{itemize}
Hence combining all the three cases we conclude $K_{\mathbb{P}(E)} + tD$ is ample for $t\geq r+2.$
\end{proof}
\end{theorem}
\begin{theorem}\label{thm3.3}
Let $\pi : \mathbb{P}(E)\longrightarrow X$ be a projective bundle corresponding to a strictly nef  vector bundle $E$ of rank $r\geq 2$ on a smooth projective variety $X$ of dimension $n$. Further assume that one of the following holds :
\begin{enumerate}
   \item $K_X+\det(E)$ is strictly nef.
  \item $E$ is ample and $K_X+\det(E) \equiv 0.$
\end{enumerate}
  Let $D\equiv m\xi+\pi^*N$ be a strictly nef divisor  on $\mathbb{P}(E)$ for some $m>0$ and some nef divisor class $N$ on $X$. Then Serrano's conjecture is true for the strictly nef divisor $D$ on $\mathbb{P}(E)$.
\begin{proof}
 We have $$K_{\mathbb{P}(E)} + tD = \bigl(tm-r\bigr)\xi+\pi^*\bigl(K_X+\det(E)+tN\bigr).$$
 Note that for $t\geq n+r+1$, we always have $(tm-r)>0$.
 Therefore, to prove  Serrano's conjecture for $D$,
it is enough to prove that for $t\geq \dim\bigl(\mathbb{P}(E)\bigr)+2$ $$\bigl(tm-r\bigr)\mu_{\min}\bigl(\eta_c^*(E\vert_C)\bigr) + K_X\cdot C + \det(E)\cdot C + tN\cdot C> 0$$
for any irreducible curve $C\subseteq X$ by Corollary \ref{corl4.5}.
\vspace{2mm}

We consider the following two cases :
\begin{itemize}
    \item Case (1)  : In this case we assume that $E$ is strictly nef and $K_X+\det(E)$ is a strictly nef line bundle. Note that  $E$ is strictly nef and thus, in particular, it is nef.  Also, by the hypothesis, we have $K_X+\det(E)$ is strictly nef. Hence $\mu_{\min}\bigl(\eta_c^*(E\vert_C)\bigr) \geq 0$ and $\bigl(K_X+\det(E)\bigr)\cdot C > 0$ for any irreducible curve $C\subseteq X$. Also $N\cdot C\geq 0$ for any irreducible $C\subseteq X$ as $N$ is nef.
    \item Case (2) : In this case we assume that $E$ is ample and $K_X+\det(E) \equiv 0$. Hence we have  $\mu_{\min}\bigl(\eta_c^*(E\vert_C)\bigr) > 0$ and $\bigl(K_X+\det(E)\bigr)\cdot C = 0$ for any irreducible curve $C\subseteq X$. Also $N\cdot C\geq 0$ for any irreducible $C\subseteq X$ as $N$ is nef.
\end{itemize}
Combining both the cases, we get that for any irreducible curve $C\subseteq X$
 \begin{align*}
  \bigl(tm-r\bigr)\mu_{\min}\bigl(\eta_c^*(E\vert_C)\bigr) + K_X\cdot C + \det(E)\cdot C +tN\cdot C> 0
 \end{align*}
This completes the proof.
\end{proof}
\end{theorem}

We now restrict our attention to projective bundles $\mathbb{P}(E)$ corresponding to special strictly nef vector bundles $E$ on smooth projective varieties $X$.

\begin{exm}
Let $X$ be a smooth complex projective variety such that $K_X$ is nef and $E=L_1\oplus L_2\oplus L_3\oplus \cdots\oplus L_r$ be a completely decomposable strictly nef vector bundle on $X$. Consider the strictly nef divisor $D=m\xi+\pi^*N$ on $\mathbb{P}(E)$ for some positive integer $m>0$ and for some nef class $N$ on $X$.

As $E$ is strictly  nef, by \cite[Proposition 2.2 (3)]{LOY19} each $L_i$ is a strictly nef line bundle on $X$. Thus
$$\det(E) = \bigotimes\limits_{i=1}^rL_i$$ is strictly nef.  Then in this case we have $K_X+\det(E)$ is strictly nef. Thus Serrano's conjecture is true for the divisor $D$.
\end{exm}
\begin{exm}
Let $E$ be a strictly nef vector bundle of rank $r$ on a smooth projective variety $X$. Consider  the following two cases :
\begin{itemize}
 \item Case (1) : Let $K_X$ is ample.
 \item Case (2) : Let $K_X \equiv 0$ and $\det(E)$ is a strictly nef line bundle.
\end{itemize}
In both these cases, we have $K_X+\det(E)$ is strictly nef. Therefore Serrano's conjecture is true for the strictly nef divisor $D=m\xi+\pi^*N$ on $\mathbb{P}(E)$ for some positive integer $m>0$ and for some nef class $N$ on $X$.
\end{exm}
\begin{exm}
\rm Let $X$ be an abelian variety so that $\
\omega_X\cong \mathcal{O}_X$ (i.e. $K_X=0$). A vector bundle $E$ on an abelian variety $X$ is called weakly-translation invariant (semi-homogeneous in the sense of Mukai) if
 for every closed point $x \in X$, there is a line bundle $L_x$ on $X$ depending on $x$ such that $T_x^*(E) \simeq E \otimes L_x$ for all $x\in X$, where $T_x$ is the translation morphism given by $x\in X$.  A semi-homogeneous vector bundle $E$ of rank $r$ on an abelian variety $X$ is ample if and only if $\det(E)$ is ample if and only if $\det(E)$ is strictly nef (cf. \cite[Corollary 3.1] {MR21}).
 Hence if $E$ is a semi-homogeneous vector bundle $E$ on an abelian variety $X$ such that $E$ is ample, and $\pi:\mathbb{P}(E)\longrightarrow X$ is the projection, then the Serrano's conjecture is true for the strictly nef divisors $D\equiv m\xi+\pi^*N$ for some positive integer $m>0$ and for some nef class $N$ on $X$.
\end{exm}
\begin{defi}
A vector bundle $E$ on a smooth projective variety $X$ is called \textit{Fano} if $\mathbb{P}(E)$ is a Fano variety,
meaning the anti-canonical line bundle $-K_{\mathbb{P}(E)}$ is ample.
\end{defi}
There are classifications of Fano bundles of rank two on smooth surfaces and on 3-folds in \cite{SW290} and \cite{SW190}
respectively. For $n\,\ge\, 4$, \cite[Main Theorem
2.4(1)]{APW94} shows that any Fano vector bundle of rank 2 on $\mathbb{P}^n$ is split. We will now only consider some of examples of non-split Fano bundles $E$ of rank two on projective spaces $\mathbb{P}^2$ and $\mathbb{P}^3$. In each of these non-split cases, we will show that $E(1)$ is nef, but not ample, and $K_{\mathbb{P}^n}+\det(E(2))$ satisfies  the hypothesis of Theorem \ref{thm3.3} for $n=2,3$. Thus  by Theorem \ref{thm3.3} the Serrano's conjecture is true for the strictly nef divisors $D\equiv m\xi+\pi^*N$ for some positive integer $m>0$ and for some nef class $N$ on $\mathbb{P}^n$, where $\xi \equiv \mathbb{P}\bigl(E(2)\bigr)$ and $\pi : \mathbb{P}\bigl(E(2)\bigr)\longrightarrow \mathbb{P}^n$ is the projection for $n=2,3$.
\begin{exm}\label{exm3.8}
Let $E$ be the rank two vector bundle on $\mathbb{P}^2$ which fits in the exact sequence
$$0\,\longrightarrow\, \mathcal{O}_{\mathbb{P}^2}\,\longrightarrow\, E \,\longrightarrow\, \mathcal{I}_{t}
\,\longrightarrow\, 0,$$
where $\mathcal{I}_t$ denotes the ideal sheaf of a point $t\in\mathbb{P}^2$. Note that $\mathcal{I}_t(1)$ is globally generated and hence nef. Thus being the extension of two nef sheaves, $E(1)$ is also nef. Also $\mathcal{I}_t(1)$ restricts trivially to any line in $\mathbb{P}^2$ through $t$. Therefore $E(1)$ is not ample, but $E(2)$ is ample.  In this case $K_{\mathbb{P}^2}+\det(E(2)) \equiv \mathcal{O}_{\mathbb{P}^2}(1)$ is ample.
\end{exm}
\begin{exm}\label{exm3.9}
Consider the rank two vector bundle $E$ on $\mathbb{P}^2$ fitting in the exact sequence
$$0\,\longrightarrow\, \mathcal{O}_{\mathbb{P}^2}(-1)^2\,\longrightarrow \,\mathcal{O}_{\mathbb{P}^2}^{\oplus4}
\,\longrightarrow\, E(1) \,\longrightarrow\, 0.$$
The vector bundle  $E(1)$ is globally generated and hence nef. But $E(1)$ is not ample as there are lines $L$ in $\mathbb{P}^2$ such that $E(1)\big\vert_L \,=\, \mathcal{O}_{L}\oplus \mathcal{O}_{L}(2)$. In this case $K_{\mathbb{P}^2}+\det(E(2)) \equiv \mathcal{O}_{\mathbb{P}^2}(1)$ is ample.
\end{exm}

\begin{exm}\label{exm3.10}
Consider the rank two vector bundle $E$ on $\mathbb{P}^2$ defined by
the exact sequence
$$0\,\longrightarrow\, \mathcal{O}_{\mathbb{P}^2}(-2)\,\longrightarrow\,
\mathcal{O}_{\mathbb{P}^2}^{\oplus 3} \,\longrightarrow\, E(1) \,\longrightarrow\, 0.$$
The vector bundle $E(1)$ is globally generated, see \cite[Proposition 2.1, Proposition 2.6]{SW290}, and $E(2)$ is ample.  In this case $K_{\mathbb{P}^2}+\det(E(2)) \equiv \mathcal{O}_{\mathbb{P}^2}(1)$ is ample.
\end{exm}
\begin{exm}\label{exm3.11}
Consider the null correlation bundle $\mathcal{N}$ on $\mathbb{P}^3$. We have the exact sequence
$$0\,\longrightarrow \,\mathcal{N} \,\longrightarrow\, T_{\mathbb{P}^3}(-1)\,\longrightarrow\,
\mathcal{O}_{\mathbb{P}^3}(1)\,\longrightarrow\, 0.$$ The vector bundle  $\mathcal{N}(1)$ is globally generated (see \cite[Theorem 2.1]{SW190}), and $N(2)$ is ample. In this case $K_{\mathbb{P}^3}+\det(N(2)) \equiv 0$.
\end{exm}
Next we consider an ample vector bundle  of rank 2 on $\mathbb{P}^2$ which is not a Fano bundle on $\mathbb{P}^2$.
\begin{exm}
    Let $Z\,=\,\{p_1,\,\cdots,\, p_5\}$ be five distinct points in $\mathbb{P}^2$ such that no
three of them are collinear. Denote the ideal sheaf of $Z$ by $\mathcal{I}_Z$.
Then there exists a vector bundle $E$ of rank 2 on $\mathbb{P}^2$ which fits in the
following short exact sequence:
$$
0 \,\longrightarrow\, \mathcal{O}_{\mathbb{P}^2}\,\stackrel{\cdot s}{\longrightarrow}\, E
\,\longrightarrow\, \mathcal{I}_Z \,\longrightarrow\, 0.$$
In other words, there is a section $s\, \in\, H^0(\mathbb{P}^2,\, E)$ whose zero locus is
precisely $Z$. Such vector bundles exist for any finite set of points in $\mathbb{P}^2$; see
\cite[Section 5.2, Example 1, Page 103]{OSS}. Let $C$ be the unique conic passing through all the points in $Z$. We then have the exact sequence
$$ 0 \,\longrightarrow\, \mathcal{O}_{C}(F) \,\longrightarrow\, E(2)\big\vert_C \,\longrightarrow\,
\mathcal{O}_C(F') \,\longrightarrow\, 0,$$
where $\deg(F) = 9, \deg(F') =-1$.
Since $E(2)\big\vert_C$ has a quotient bundle of negative degree, it is not nef.
Consequently, $E(2)$ is not a nef line bundle on $\mathbb{P}^2$.
It can be shown that $E(3)$ is ample. In this case also, $K_{\mathbb{P}^2}+\det\bigl(E(3)\bigr) \equiv \mathcal{O}_{\mathbb{P}^2}(3)$ is ample. Thus by Theorem \ref{thm3.3} the Serrano’s conjecture is true for the strictly nef
divisors $D \equiv m\xi + \pi^*N$ for some positive integer $m > 0$ and for some nef class $N$ on $\mathbb{P}^2$, where $\xi \equiv \mathcal{O}_{\mathbb{P}(E(3))}(1)$ and $\pi : \mathbb{P}(E(3)) \longrightarrow \mathbb{P}^2$ is the projection.
\end{exm}
%\begin{defi}
%For a vector bundle $E$ of rank $r$ on a smooth projective variety $X$, the element $$c_2\bigl(\End(E)\bigr) = 2rc_2(E)-(r-1)c_1^2(E) \in H^4(X,\mathbb{Q})$$ is called the \it discriminant \rm of $E$. If $E$ is a $\mu_H$- semistable vector bundle on a smooth projective variety $X$ of dimension $n$, then Bogomolov's inequality says that $c_2\bigl(\End(E)\bigr)\cdot H^{n-2}\geq 0$.
%\end{defi}
\begin{theorem}\label{thm4.9}
 Let $E$ be a semistable vector bundle or rank $r$ on a smooth projective irreducible variety $X$ with $c_2(\End(E))\equiv 0$. Then $E$ is striclty nef vector bundle if $\det(E)$ is strictly nef.
 \begin{proof}
Let $\phi:C\longrightarrow X$ be a non-constant finite morphism from a smooth curve $C$ to $X$, and $$\phi^*(E)\longrightarrow L\longrightarrow 0$$ be a quotient line bundle of $\phi^*(E)$. As $E$ is a semistable vector bundle with $c_2\bigl(\End(E)\bigr)=0$, by \cite[Theorem 3.10]{FL22} we have $\phi^*(E)$ is semistable. Thus by definition of semistability $$\deg(L)=\mu(L) \geq \mu(\phi^*(E)).$$ As $\det(E)$ is strictly nef, we have $\det(E)\cdot \phi(C) >0$. Thus $$\deg(L)=\mu(L) \geq \mu(\phi^*(E)) = \frac{\det(E)\cdot \phi(C)}{\rk(E)} > 0.$$
Hence by \cite[Proposition 2.1]{LOY19}, the semistable vector bundle $E$ is striclty nef.
 \end{proof}
\end{theorem}
\begin{xrem}
  By Mumford's example \cite[Example 2.5, Example 2.7]{LOY21}, there are strictly nef vector bundle $E$ on smooth projective curves $C$ such that $\deg(E) =0$. Since $\det(E)$ is numerically trivial, it is not strictly nef.
\end{xrem}
\begin{theorem}\label{thm4.10}
 Let $E$ be a semistable vector bundle or rank $r$ on a smooth projective irreducible variety $X$ of dimension $n$ with $c_2(\End(E))=0$, and  $\det(E)$ is strictly nef. Consider the projection map $\pi:\mathbb{P}(E)\longrightarrow X$, and the strictly nef divisor $D\equiv m\xi+\pi^*N$ for some integer $m\geq r$ and for some nef class $N$ on $X$. Then Serrano's conjecture is true for $D$.
 \begin{proof}
 By Corollary \ref{corl4.5}
 \begin{align*}
 K_{\mathbb{P}(E)}+tD = \bigl(tm-r\bigr)\xi+\pi^*(K_X+\det(E)+tN)
\end{align*}
is ample if and only if $\bigl(tm-r\bigr) >0$ and $$\bigl(tm-r\bigr)\mu(\bigl(\eta_c^*(E\vert_C)\bigr) + K_X\cdot C+\det(E)\cdot C+tN\cdot C >0.$$ for every irreducible curve $C\subseteq X$. Now for $t\geq n+r+1$, we have $\bigl(tm-r\bigr) >0$. Let $C\subseteq X$ be an irreducible curve. As $\det(E)$ is strictly nef, then by \cite[Lemma 1.1]{S95} $K_X+t\det(E)$ is strictly nef for $t\geq n+2.$ Also $\eta_c^*(E\vert_C)$ is a semistable bundle for any irreducible curve $C\subseteq X$ so that $$\mu_{\min}\bigl(\eta_c^*(E\vert_C)\bigr) = \mu\bigl(\eta_c^*(E\vert_C)\bigr) = \frac{\det(E)\cdot C}{\rk(E)} > 0,$$
as $\det(E)$ is strictly nef. Therefore
\begin{align*}
 &\bigl(tm-r\bigr)\mu(\bigl(\eta_c^*(E\vert_C)\bigr) + K_X\cdot C+\det(E)\cdot C+tN\cdot C\\
 &=\bigl(tm-r\bigr)\frac{\det(E)\cdot C}{r} + K_X\cdot C+\det(E)\cdot C+tN\cdot C\\
 &=\bigl(\frac{tm}{r}-1\bigr)\det(E)\cdot C+K_X\cdot C+\det(E)\cdot C+tN\cdot C \\
 &=\frac{tm}{r}\det(E)\cdot C+K_X\cdot C+tN\cdot C
\end{align*}
As $t\geq n+r+1$ and $m\geq r$, we have $$\frac{tm}{r}\geq n+r+1\geq n+2.$$ Thus for $t\geq n+r+1$ we get $\frac{tm}{r}\det(E)+K_X$ is strictly nef, i.e. $$\frac{tm}{r}\det(E)\cdot C+K_X\cdot C>0.$$
Also $N$ being nef, we have $N\cdot C \geq 0$.
This shows that $$\bigl(tm-r\bigr)\mu(\bigl(\eta_c^*(E\vert_C)\bigr) + K_X\cdot C+\det(E)\cdot C+tN\cdot C > 0$$ for every irreducible curve $C$ in $X$ for $t\geq n+r+1$. Therefore $K_{\mathbb{P}(E)}+tD$ is ample for $t\geq n+r+1$. This completes the proof.
 \end{proof}
\end{theorem}
\begin{theorem}\label{thm4.11}
Let $E$ be a semistable vector bundle or rank $r$ on a smooth projective irreducible variety $X$ with $c_2(\End(E))=0$. Further assume that $-K_X$ is nef. Then the following are equivalent :
\begin{itemize}
 \item (i) $E$ is ample.
 \item (ii) $\det(E)$ is ample.
 \item (iii) $\det(E)$ is strictly nef.
\end{itemize}
\begin{proof}
$(i) \iff (ii)$ This equivalence follows from \cite[Theorem 1] {MR21}.

$(ii) \implies (iii)$ This follows from the definition of strictly nef line bundle and Nakai-Moishezon-Kliemann criterion for ampleness.

$(iii) \implies (i)$ Let $\det(E)$ is strictly nef on $X$ and $\pi:\mathbb{P}(E)\longrightarrow X$ be the projection. By Theorem \ref{thm4.10}, we have $K_{\mathbb{P}(E)}\otimes L^{\otimes t}$ is ample for $t$ large enough integer, where $L=\mathcal{O}_{\mathbb{P}(E)}(1)$. Since $-K_X$ is nef, we conclude that $K_{\mathbb{P}(E)/X}\otimes L^{\otimes t}$ is ample for large enough $t$.
Thus $\det\bigl(\pi_*(K_{\mathbb{P}(E)/X}\otimes(K_{\mathbb{P}(E)/X}\otimes L^{\otimes t}))\bigr)$ is ample, and hence $\det(E)$ is ample.
\end{proof}
\end{theorem}
\begin{theorem}\label{thm4.12}
Let $\pi: Y= \mathbb{P}(V)\longrightarrow X$ be a projective bundle corresponding to a vector bundle $E$ of rank $r$ on a smooth projective  variety $X$. Then $-K_{\mathbb{P}(V)}$ is nef if and only if $$0\leq \det(V)\cdot C-r\mu_{\min}\bigl(\eta_c^*(V\vert_C)\bigr)\leq -K_X\cdot C$$ for every irreducible curve $C\subseteq X$, where $\eta_c:\tilde{C}\longrightarrow  C$ is the normalization map. In particular, if $K_X\equiv 0$, then $-K_{\mathbb{P}(V)}$ is nef if and only if $\eta_c^*(V\vert_C)$ is semistable for every irreducible curve $C\subseteq X$.
\begin{proof}
We note that $$-K_{\mathbb{P}(V)} = r\xi + \pi^*\bigl(-K_X-\det(V)\bigr).$$ Hence by Remark \ref{xrem1} $-K_{\mathbb{P}(V)}$ is nef if and only if $$r\mu_{\min}\bigl(\eta^*_c(V\vert_C)\bigr) -K_X\cdot C - \det(V)\cdot C \geq 0$$ for every irreducible curve $C\subseteq X$. Note that for any irreducible curve $C\subseteq X$ we have $$\det(V)\cdot C = \deg\bigl(\eta_c^*(V\vert_C)\bigr) \geq r\mu_{\min}\bigl(\eta^*_c(V\vert_C)\bigr).$$
Therefore, $-K_{\mathbb{P}(V)}$ is nef if and only if $$0\leq \det(V)\cdot C-r\mu_{\min}\bigl(\eta_c^*(V\vert_C)\bigr)\leq -K_X\cdot C$$ for every irreducible curve $C\subseteq X$.
This implies that if $K_X\cdot C >0$ for some irreducible curve $C\subseteq X$, then $-K_{\mathbb{P}(V)}$ is never nef.

Now, in particular, if $K_X\equiv 0$, then $-K_X\cdot C = 0$ for every irreducible curve $C\subseteq X$. Thus in this case, $-K_{\mathbb{P}(V)}$ is nef if and only if for every irreducible curve $C\subseteq X$ $$\det(V)\cdot C = r\mu_{\min}\bigl(\eta_c^*(V\vert_C)\bigr),$$ or equivalently $\mu\bigl(\eta_c^*(V\vert_C)\bigr) = \mu_{\min}\bigl(\eta_c^*(V\vert_C)\bigr),$ i.e.  $\eta_c^*(V\vert_C)$ is semistable for every irreducible curve $C\subseteq X$.
\end{proof}
\end{theorem}
\begin{xrem}
\rm Let $\pi: X=\mathbb{P}(V)\longrightarrow C$ be a projective bundle over a smooth irreducible projective curve $C$. Let us fix the notations : $\xi\equiv \mathcal{O}_{\mathbb{P}(V)}(1)$, $f\equiv$  a fiber of $\pi$. Then by Theorem \ref{thm4.12} $-K_{\mathbb{P}(V)}$ is nef if and only if $$0\leq \deg(V)-r\mu_{\min}(V) \leq 2-2g,$$ where $r$ is the rank of the vector bundle $V$ and $g$ is the genus of $X$.
Hence $-K_{\mathbb{P}(V)}$ is not nef if $g\geq 2.$
\vspace{1mm}

For $g=1$, i.e when $C$ is an elliptic curve, the anti-canonical bundle $-K_{\mathbb{P}(V)}$ is nef if and only if $V$ is a semistable bundle on $C$. Let $V$ be a semistable bundle of rank $r$ on an elliptic curve $C$, and $E$ is a semistable vector bundle of rank $l$ on $\mathbb{P}(V)$ with $c_2\bigl(\End(E)\bigr)=0$.
Then by Theorem \ref{thm4.11} and \cite[Corollary 7]{MR21} the line bundle $\det(E) \equiv x\xi+yf\in N^1(X)_{\mathbb{R}}$, is strictly nef if and only if $\det(E)$ is ample if and only if $x>0$ and $x\mu(V)+y>0$.
\vspace{2mm}

Now consider the case $C=\mathbb{P}^1$. Let $V$ be a vector bundle or rank $r$ on $\mathbb{P}^1$ such that $$\deg(V)-r\mu_{\min}(V)\leq 2.$$
 Then $-K_{\mathbb{P}(V)}$ is nef.  Let $E$ be a semistable vector bundle on $\mathbb{P}(V)$ with $c_2\bigl(\End(E)\bigr)=0$. Then by Theorem \ref{thm4.11} and \cite[Corollary 7]{MR21} the determinant line bundle $\det(E)\equiv x\xi+yf\in N^1(X)_{\mathbb{R}}$ is strictly nef if and only if $\det(E)$ is ample if and only if $x>0$ and $x\mu_{\min}(V)+y>0$.
\end{xrem}

\begin{corl}
Let $V$ be a semi-homogeneous vector bundle of rank $r$ on an abelian variety $X$ and $\pi:Y:= \mathbb{P}(V)\longrightarrow X$ be the projection. Let $E$ be a semistable vector bundle on $Y$ with $c_2\bigl(\End(E)\bigr)=0$. We fix the notaions $\xi\equiv \mathcal{O}_{\mathbb{P}(V)}(1)$. Then $\det(E)\equiv x\xi+\pi^*L$ for some $L\in N^1(X)_{\mathbb{R}}$ is strictly nef if and only if $x>0$ and $\frac{x}{r}\det(V)+L$ is strictly nef on $X$.
\begin{proof}
A result due to Mukai proved that a semi-homogeneous vector bundle $V$ on an abelian variety $X$ is Gieseker semistable (see \cite[Chapter 1]{H-L} for definition) with respect to some polarization and it has projective Chern classes zero, i.e., if $c(V)$ is the total Chern class, then $c(V) = \Bigl\{ 1+ c_1(V)/r\Bigr\}^r$ (see \cite[Theorem 5.8, page 260]{Muk78}, \cite[ Proposition 6.13, page 266]{Muk78}. Gieseker semistablity implies slope semistability (see \cite{H-L}). So, in particular, we have $V$ is slope semistable with $c_2\bigl(\End(V)\bigr) = 2rc_2(V) - (r -1)c_1^2(V) = 0$. Therefore $\eta_c^*(V\vert_C)$ is semistable for any irreducible curve $C\subseteq X$, where $\eta_c : \tilde{C}\longrightarrow C$ is the normalization map. Hence $$\mu_{\min}\bigl(\eta_c^*(V\vert_C)\bigr) = \mu\bigl(\eta_c^*(V\vert_C)\bigr) = \frac{\det(E)\cdot C}{r}$$
for every irreducible curve $C$ in $X$. Also, by Theorem \ref{thm4.12} we have $-K_{\mathbb{P}(V)}$ is nef as $K_X =0$.
Therefore, by Theorem \ref{thm4.11}, the line bundle $\det(E) \equiv x\xi+\pi^*L$ is strictly nef if and only if $\det(E)$ is ample if and only if $x>0$ and $x\mu_{\min}\bigl(\eta_c^*(V\vert_C)\bigr)+L\cdot C>0$ for every irreducible curve $C\subseteq X$.
This is equivalent to the fact that $x>0$ and $\frac{x}{r}\det(V)+L$ is strictly nef on $X$.
\end{proof}
\end{corl}
\section{Results on Toric varieties}
Let $X$ be a toric variety, and $E$ be a torus equivariant vector bundle of rank $r$ on $X$. Let $l_1,l_2,\cdots,l_m$ be finitely many invariant curves in $X$. Then  each invariant curve $l_i$ is rational, i.e. $l_i\cong \mathbb{P}^1$ for each $i\in\{1,2,\cdots,m\}$. For a fixed invariant curve $l_i$, consider the following commutative fiber product  diagram:
\begin{center}
 \begin{tikzcd}
 \mathbb{P}(E)\times_Xl_i = \mathbb{P}(E\vert_{l_i}) \arrow[r, "\psi_i"] \arrow[d,"\pi_i"] & \mathbb{P}(E) \arrow[d,"\pi"] \\
         l_i \arrow[r, "i"]                                         & X
\end{tikzcd}
\end{center}
where $i:l_i\hookrightarrow X$ is the inclusion. Let $\xi_i$ and $f_i$ be the numerical classes of the tautological line bundle $\mathcal{O}_{\mathbb{P}(E\vert_{l_i})}(1)$ and a fiber of the map $\pi_i$ respectively. Then
\begin{align*}
 \overline{\NE}\bigl(\mathbb{P}(E\vert_{l_i})\bigr) = \Bigl\{a(\xi_i^{r-2}f_i)+b(\xi^{r-1}-m_i\xi^{r-2}f_i)\mid a,b\in \mathbb{R}_{\geq 0}\Bigr\},
\end{align*}
where $m_i:= \deg(E\vert_{l_i})-\mu_{\min}(E\vert_{l_i})$. For each $i\in \{1,2,\cdots,m\},$ let us fix the following notations :
\vspace{2mm}

\begin{center}
 $\Sigma_i:= \xi_i^{r-2}f_i\hspace{1mm}, \hspace{4mm} \Omega_i:= \xi^{r-1}-m_i\xi^{r-2}f_i$
 \end{center}
 and
 \begin{center}
 $C_0:= \psi_i(\Sigma_i)\hspace{1mm}, C_i := \psi_i(\Omega_i)$ for every $i=1,2,\cdots,m$.
 \end{center}
 \vspace{2mm}

 Then we have the following result (see \cite[Proposition 3.2]{DKS22}).
 \begin{prop}\label{corl5.1}
  With notations as above, the closed cone of curves
  $$\overline{\NE}\bigl(\mathbb{P}(E)\bigr) = \Bigl\{a_0C_0+a_1C_1+a_2C_2+\cdots+a_mC_m\mid a_i\in \mathbb{R}_{\geq 0}\Bigr\}.$$
\end{prop}
We use the above Proposition \ref{corl5.1} to find the ample cone $\Amp\bigl(\mathbb{P}(E)\bigr)$ of $\mathbb{P}(E)$.
\begin{prop}\label{prop5.2}
Let $E$ be a torus equivariant vector bundle of rank $r$  on a toric variety $X$, and $\pi:\mathbb{P}(E)\longrightarrow X$ be the projection. Then an element $y_0\xi+\pi^*(\gamma)$ is in the ample cone $\Amp\bigl(\mathbb{P}(E)\bigr)$ if and only if $y_0>0$ and $y_0\mu_{\min}(E\vert_{l_i})+\gamma\cdot l_i>0$ for every $i\in\{1,2,\cdots,m\}$.
\begin{proof}
   The proof is similar to the proof of Corollary \ref{corl4.5}. We omit the details.
\end{proof}
\end{prop}
\begin{theorem}
 Let $E$ be a torus equivariant strictly nef vector bundle of rank $r$  on a toric variety $X$ of dimension $n$, and $\pi:\mathbb{P}(E)\longrightarrow X$ be the projection. Further assume that $K_X$ is nef. Consider the strictly nef divisor $D\equiv m\xi+\pi^*N$ for some positive integer $m>0$ and for some nef class $N$ on $X$. Then Serrano's conjecture is true for $D$.
 \begin{proof}
As $E$ is strictly nef,  by \cite[Proposition 2.2(1)]{LOY19} it's restriction $E\vert_{l_i}$ is strictly nef for every invariant curve $l_i\subseteq X$ in $X$. Thus as $l_i\cong \mathbb{P}^1$, we have $E\vert_{l_i}$ is ample.
Hence $\mu_{\min}(E\vert_{l_i})>0$ for every invariant curve $l_i$ in $X$. Also, we have $K_X\cdot l_{i} \geq 0$ and $N\cdot l_{i}\geq 0$ for every invariant curve $l_i\subseteq X$  a.s $K_X$ and $N$ are nef.

Now
$$K_{\mathbb{P}(E)}+tD = \bigl(tm-r\bigr)\xi+\pi^*(K_X+\det(E)+tN)$$
is ample if and only if $\bigl(tm-r\bigr)>0$ and for every invariant curve $l_i\subseteq X$
$$\bigl(tm-r\bigr)\mu_{\min}(E\vert_{l_i})+K_X\cdot l_i+\det(E)\cdot l_i+tN\cdot l_i > 0$$
For $t\geq r+n+1$, we have $\bigl(tm-r\bigr)>0$. Also for any invariant curve $l_i\subseteq X$
\begin{align*}
 &\bigl(tm-r\bigr)\mu_{\min}(E\vert_{l_i})+K_X\cdot l_i+\det(E)\cdot l_i+tN\cdot l_i \\
 &= tm\mu_{\min}(E\vert_{l_i})+\deg(E\vert_{l_i})-r\mu_{\min}(E\vert_{l_i})+K_X\cdot l_i+tN\cdot l_i >0.
\end{align*}
Therefore we conclude that Serrano's conjecture is true for the strictly nef divisor $D$.
\end{proof}
\end{theorem}
Next we restrict our attention to projective bundle $\mathbb{P}(E)$ corresponding to torus equivariant strictly nef vector bundle $E$ on projective spaces $\mathbb{P}^n$. There are $\binom{n+1}{2}$ many invariant curves in $\mathbb{P}^n$ for $n\geq 2$. Each invariant curve in $\mathbb{P}^n$ is a line. Since the Picard number of $\mathbb{P}^n$ is $1$, we have any two invariant curves in $\mathbb{P}^n$ are numerically equivalent.

\begin{theorem}\label{thm4.4}
Let $E= \mathcal{O}_{\mathbb{P}^n}(a_1)\oplus \mathcal{O}_{\mathbb{P}^n}(a_2)\oplus \mathcal{O}_{\mathbb{P}^n}(a_3)\oplus\cdots\oplus \mathcal{O}_{\mathbb{P}^n}(a_r)$ be a completely decomposable strictly nef vector bundle of rank $r\geq 2$ on $\mathbb{P}^n$ such that $ a_1\leq a_2 \leq a_3\leq \cdots\leq a_r$. Consider the strictly nef divisor $D=m\xi+\pi^*N$ for some positive integer $m>0$ and some nef class $N$ on $\mathbb{P}^n$. Then Serrano's conjecture is true for $D$.
\begin{proof}
We denote the numerical class of $\pi^*\bigl(\mathcal{O}_{\mathbb{P}^n}(1)\bigr)$ by $H$. Let $N\equiv \mathcal{O}_{\mathbb{P}^n}(l)$ for some $l\geq 0$.
Then
\begin{align*}
& K_{\mathbb{P}(E)}+tD = -r\xi+\pi^*\bigl(K_{\mathbb{P}^n}+\det(E)\bigr) +t(m\xi+\pi^*N)\\
&= (tm-r)\xi+\bigl(\sum\limits_{i=1}^ra_i-(n+1)+l\bigr)H
\end{align*}
Now for $t\geq r+n+1,$ we have $(tm-r) >0$ as $m>0$.
By \cite[Proposition 2.2 (3)]{LOY19} each direct summand $\mathcal{O}_{\mathbb{P}^n}(a_i)$ is strictly nef and hence $a_i>0$ for each $i\in\{1,2,\cdots,r\}$. Thus for $t\geq r+n+1$ we get
\begin{align*}
 & (tm-r)a_1+\bigl(\sum\limits_{i=1}^ra_i-(n+1)+l\bigr)\\
&= a_1tm-a_1r +a_1+ a_2+a_3+a_3+\cdots+a_r-(n+1)+l\\
&\geq a_1(r+n+1) -a_1r + a_2 + a_3+\cdots+a_r - (n+1)+l\\
&=\bigl((n+1)a_1-(n+1)\bigr)+a_2+a_3+\cdots+a_r+l\\
&= (n+1)(a_1-1)+a_2+a_3+\cdots+a_r+l\\
&>0.
\end{align*}
Thus by Proposition \ref{prop5.2} $K_{\mathbb{P}(E)}+tD$ is ample for $t\geq r+n+1$.
 \end{proof}
\end{theorem}
\begin{theorem}\label{thm4.5}
Let $\pi:\mathbb{P}(T_{\mathbb{P}^n})\longrightarrow \mathbb{P}^n$ be the projection, where $T_{\mathbb{P}^n}$ is the tangent bundle on $\mathbb{P}^n$. Consider the strictly nef divisor $D\equiv m\xi+\pi^*N$ for some positive integer $m>0$ and for some nef class $N$ on $\mathbb{P}^n$. Then Serrano's conjecture is true for $D$.
\begin{proof}
We note that  $T_{\mathbb{P}^n}$ is ample and hence it is strictly nef. Also its determinant line bundle $\det(T_{\mathbb{P}^n})$ is an ample line bundle and thus $\det(T_{\mathbb{P}^n})\cdot l_i>0$ for every invariant curve $l_i$.

For any invariant curve $l_i\subset \mathbb{P}^n$, we have $T_{\mathbb{P}^n}\vert_{l_i}\cong \mathcal{O}_{\mathbb{P}^1}(2)\oplus \mathcal{O}_{\mathbb{P}^1}(1)\oplus \cdots \oplus \mathcal{O}_{\mathbb{P}^1}(1).$
Hence $\mu_{\min}(T_{\mathbb{P}^n}\vert_{l_i}) = 1$ for every invariant curve $l_i\subseteq \mathbb{P}^n$.
\vspace{1mm}

For $t\geq \dim\bigl(\mathbb{P}(T_{\mathbb{P}^n})\bigr) + 2 = 2n+1,$ we have $(tm-n) >0$, and for every invariant curve $l_i\subseteq \mathbb{P}^n$
\begin{align*}
&\bigl(tm-n)\mu_{\min}(T_{\mathbb{P}^n}\vert_{l_i})+K_{\mathbb{P}^n}\cdot l_i+\det(E)\cdot l_i +tN\cdot l_i \\
&=tm+K_{\mathbb{P}^n}\cdot l_i +\deg(E\vert_{l_i})-n\mu_{\min}(T_{\mathbb{P}^n}+tN\cdot l_i \\
& \geq 2n+1-(n+1)+\bigl(\deg(E\vert_{l_i})-n\mu_{\min}(T_{\mathbb{P}^n})\bigr)+tN\cdot l_i\\
&=n+\bigl(\deg(E\vert_{l_i})-n\mu_{\min}(T_{\mathbb{P}^n})\bigr)+tN\cdot l_i > 0.
\end{align*}
Therefore we conclude by Proposition \ref{prop5.2} that $K_{\mathbb{P}(\mathbb{T}_{\mathbb{P}^n})}+tD$ is ample for $t\geq 2n+1$. This completes the proof.
\end{proof}
\end{theorem}
In view of Theorem \ref{thm4.4}, Theorem \ref{thm4.5} and Example \ref{exm3.8} - Example \ref{exm3.11}, we conclude the following.
\begin{theorem}
Let $E$ be a Fano ample vector bundle of rank 2 on $\mathbb{P}^n$ (which means $\mathbb{P}(E)$ is a Fano variety), and $\pi : \mathbb{P}(E)\longrightarrow \mathbb{P}^n$ be the projection. Let $D \equiv m\xi + \pi^*N$ be a strictly nef divisor on $\mathbb{P}(E)$ for some $m > 0$ and some nef divisor class $N$ on $\mathbb{P}^n$.
Then Serrano’s conjecture is true for the strictly nef divisor $D$ on $\mathbb{P}(E)$.
\end{theorem}
\section{Acknowledgement}
 The  author is supported financially by SERB-National Postdoctoral Fellowship (File no : PDF/2021/00028).

%\bibliographystyle{amsplain}


\end{document}